\documentclass{article}

\title{\LARGE \textbf{Disconnected Forbidden Subgraphs, Toughness and Hamilton Cycles}}
\author{Zh.G. Nikoghosyan\footnote{G.G. Nicoghossian (up to 1997)}}

\begin{document}

\maketitle

\begin{abstract}
In 1974, Goodman and Hedetniemi proved that every 2-connected $(K_{1,3},K_{1,3}+e)$-free graph is hamiltonian. This result gave rise many other hamiltonicity conditions  for various pairs and triples of forbidden connected subgraphs under additional connectivity conditions.  In 1997, it was proved that  a single forbidden connected subgraph $R$ in 2-connected graphs can create only a trivial class of hamiltonian graphs (complete graphs) with $R=P_3$. In this paper we prove that a single forbidden subgraph $R$ can create a non trivial class of hamiltonian graphs if $R$ is disconnected: $(\ast1)$ every $(K_1\cup P_2)$-free graph either is hamiltonian or belongs to a well defined class of non hamiltonian graphs; $(\ast2)$ every 1-tough $(K_1\cup P_3)$-free graph is hamiltonian.  We conjecure that  every 1-tough $(K_1\cup P_4)$-free graph is hamiltonian and every 1-tough $P_4$-free graph is hamiltonian\\


\end{abstract}

\section{Introduction}

Only finite undirected graphs without loops or multiple edges are considered. We denote by $n$ and $\alpha$ the order and the independent number of a graph, respectively. If $H_1,...,H_t$ $(t\ge1)$ are graphs then a graph $G$ is said to be $(H_1,...,H_t)$-free if $G$ contains no copy of any of the graphs $H_1,...,H_t$. Further, we denote by $P_i$ and $C_i$ the path and the cycle on $i$ vertices.  A good reference for any undefined terms is \cite{[1]}. 

The first sufficient condition for hamiltonicity of a graph in terms of forbidden subgraphs is due to Goodman and Hedetniemi \cite{[3]}.\\

\noindent\textbf{Theorem A}. Every 2-connected $(K_{1,3}, K_{1,3}+e)$-free graph is hamiltonian.\\

This result gave rise many other hamiltonicity conditions  for various pairs and triples of forbidden connected subgraphs under additional connectivity conditions. 

In 1997, Faudree and Gould \cite{[2]} proved that a single forbidden connected subgraph $R$ in 2-connected graphs can create only a trivial class (complete graphs) of hamiltonian graphs with $R=P_3$. \\

\noindent\textbf{Theorem B}. Let $R$  be a connected graph  and $G$ be a 2-connected graph. Then $G$ is $R$-free implies $G$ is hamiltonian if and only if $R=P_3$. \\

In this paper we prove that a single forbidden subgraph $R$ can create a non trivial class of hamiltonian graphs if $R$ is disconnected. First of all, observe the following.\\

\noindent\textbf{Proposition 1}. Every $(K_1\cup K_1)$-free graph is complete and therefore is hamiltonian. \\

It is not hard to see that every $(K_1\cup K_1\cup K_1)$-free graph either is hamiltonian or consists of two complete graphs having at most one vertex in common. In other words, we have the following.\\

\noindent\textbf{Proposition 2}. Every 2-connected $(K_1\cup K_1\cup K_1)$-free graph is hamiltonian.\\

Observe that $K_1\cup P_2$ is the minimum forbidden disconnected subgraph containing at least one edge. To describe the hamiltonian graphs in  $(K_1\cup P_2)$-free graphs, we need the following recursive definition.\\

\noindent\textbf{Definition}. We say that $G\in \aleph$ if and only if either $V(G)$ is independent set of vertices or $G$ is complete graph or there is a bipartition $V=V_1\cup V_2$ such that

1) $V_1$ is an independent set of vertices,

2) $N(v)=V_2$ for each $v\in V_1$,

3) $G[V_2]\in \aleph$. \\

\noindent\textbf{Theorem 1}. Every $(K_1\cup P_2)$-free graph $G$ either is hamiltonian or $G\in \aleph$ with $\alpha(G)>n/2$.\\

The following corollary follows immediately.\\

\noindent\textbf{Corollary 1}. Every 1-tough $(K_1\cup P_2)$-free graph is hamiltonian.\\

Further relaxing of the condition "$G$ is $(K_1\cup P_2)$-free" implies the following.\\

\noindent\textbf{Theorem 2}. Every 1-tough $(K_1\cup P_3)$-free graph is hamiltonian.\\

Examples for sharpness. Clearly, $K_{2,3}$ is a non hamiltonian $(K_1\cup P_3)$-free (even $(K_1\cup P_2)$-free) graph with $\tau(K_{2,3})=2/3$ and $\kappa (K_{2,3})=2$,  implying that the condition "$G$ is 1-tough" in Theorem 2 can not be removed or replaced by "$G$ is 2-connected".  Now form a graph, denoted by $H_1$, by adding a new vertex $x_7$ to $C_6=x_1x_2x_3x_4x_5x_6x_1$ and new edges $x_7x_1,x_7x_4,x_2x_6$. Since $H_2$ is a non hamiltonian $(K_2\cup P_3)$-free graph with $\tau(H_2)=1$,  we can claim that the condition "$G$ is $(K_1\cup P_3)$-free" in Theorem 2 cannot be relaxed to  "$G$ is $(K_2\cup P_3)$-free". Finally, $H_2$ is a $(K_1\cup K_{1,3})$-free graph and hence the condition "$G$ is $(K_1\cup P_3)$-free" in Theorem 2 cannot be relaxed to  "$G$ is $(K_1\cup K_{1,3})$-free". So, Theorem 2 is best possible in many respects. The condition  "$G$ is $(K_1\cup P_3)$-free" in Theorem 2 perhaps can be relaxed to "$G$ is $(K_1\cup P_4)$-free".   \\

\noindent\textbf{Conjecture 1}. Every 1-tough $(K_1\cup P_4)$-free graph is hamiltonian.\\

Example for sharpness. The graph $H_1$ (see the sharpness examples for Theorem 2) shows that the condition "$G$ is $(K_1 \cup P_4)$-free" in Conjecture 1 (if true) can not be replaced by "$G$ is $(K_1 \cup P_5)$-free".\\

For more than one tough graphs, the following is reasonable.\\

\noindent\textbf{Conjecture 2}. Every $(K_1\cup P_5)$-free graph with $\tau>1$ is hamiltonian.\\

Example for sharpness. The Petersen graph shows that the condition "$G$ is $(K_1 \cup P_5)$-free" in Conjecture 2 can not be replaced by "$G$ is $(K_1 \cup P_6)$-free".\\

The next conjecture concerns $(K_2\cup K_2)$-free graphs.\\

\noindent\textbf{Conjecture 3}. Every $(K_2\cup K_2)$-free graph with $\tau>1$ is hamiltonian.\\

Examples for sharpness. The graph $H_1$ (see the sharpness examples for Theorem 2) shows that the condition $\tau>1$ in Conjecture 3 can not be replaced by $\tau=1$. Further, the Petersen graph shows that the condition "$G$ is $(K_2\cup K_2)$-free" in Conjecture 3 can not be replaced by "$G$ is $(K_2\cup K_3)$-free".

The next conjecture is based on $K_{1,3}$ (Claw).\\

\noindent\textbf{Conjecture 4}. Every $(K_1\cup K_{1,3})$-free graph with $\tau>4/3$ is hamiltonian.\\

Example for sharpness. The Petersen graph shows that the condition $\tau>4/3$ in Conjecture 4 can not be replaced by $\tau=4/3$.\\

Finally, we hope that the following is true.\\

\noindent\textbf{Conjecture 5}. Every 1-tough $P_4$-free graph is hamiltonian.\\

\section{Notations and preliminaries}

The set of vertices of a graph $G$ is denoted by $V(G)$ and the set of edges by $E(G)$. For $S$ a subset of $V(G)$, we denote by $G\backslash S$ the maximum subgraph of $G$ with vertex set $V(G)\backslash S$. We write $G[S]$ for the subgraph of $G$ induced by $S$. For a subgraph $H$ of $G$ we use $G\backslash H$ short for $G\backslash V(H)$. The neighborhood of a vertex $x\in V(G)$ will be denoted by $N(x)$.  Furthermore, for a subgraph $H$ of $G$ and $x\in V(G)$, we define $N_H(x)=N(x)\cap V(H)$. Let $s(G)$ denote the number of components of a graph $G$. A graph $G$ is $t$-tough if $|S|\ge ts(G\backslash S)$ for every subset $S$ of the vertex set $V(G)$ with $s(G\backslash S)>1$. The toughness of $G$, denoted $\tau(G)$, is the maximum value of $t$ for which $G$ is $t$-tough (taking $\tau(K_n)=\infty$ for all $n\ge 1$).

A simple cycle (or just a cycle) $C$ of length $t$ is a sequence $v_1v_2...v_tv_1$ of distinct vertices $v_1,...,v_t$ with $v_iv_{i+1}\in E(G)$ for each $i\in \{1,...,t\}$, where $v_{t+1}=v_1$. When $t=2$, the cycle $C=v_1v_2v_1$ on two vertices $v_1, v_2$ coincides with the edge $v_1v_2$, and when $t=1$, the cycle $C=v_1$ coincides with the vertex $v_1$. So, all vertices and edges in a graph can be considered as cycles of lengths 1 and 2, respectively. A graph $G$ is hamiltonian if $G$ contains a Hamilton cycle, i.e. a cycle of length $n$. A cycle $C$ in $G$ is dominating if  $G\backslash C$ is edgeless.

Paths and cycles in a graph $G$ are considered as subgraphs of $G$. If $Q$ is a path or a cycle, then the length of $Q$, denoted by $|Q|$, is $|E(Q)|$. We write $Q$ with a given orientation by $\overrightarrow{Q}$. For $x,y\in V(Q)$, we denote by $x\overrightarrow{Q}y$ the subpath of $Q$ in the chosen direction from $x$ to $y$. For $x\in V(C)$, we denote the $h$-th successor and the $h$-th predecessor of $x$ on $\overrightarrow{C}$ by $x^{+h}$ and $x^{-h}$, respectively. We abbreviate $x^{+1}$ by $x^+$. For each $X\subset V(C)$, we define $X^{+h}=\{x^{+h}|x\in X\}$. \\

\section{Proofs}

\noindent \textbf{Proof of Theorem 1}. Let $G$ be a $(K_1\cup P_2)$-free graph. If $V(G)$ is independent then $G\in \aleph$ and we are done. Let $G$ contains at least one edge. Next, if $G$ is disconnected then clearly $G$ contains $K_1\cup P_2$ as in induced subgraph, contradicting the hypothesis. Let $G$ is connected. Further, if $G$ is a tree then clearly $G$ either is a star (that is $G$ is a complete bipartite graph and hence $G\in\aleph$) or contains $K_1\cup P_2$ as in induced subgraph, again contradicting the hypothesis. Now let $G$ is not a tree, that is contains a cycle, and let $C=v_1v_2...v_tv_1$ be a longest cycle in $G$. If $V(G\backslash C)=\emptyset$ then $C$ is a Hamilton cycle and we are done. Let $V(G\backslash C)\not=\emptyset$. It follows that $xy\in E(G)$ for some $x\in V(G\backslash C)$ and $y\in V(C)$.  Assume w.l.o.g. that $y=v_1$. Since $C$ is extreme, we have $xv_2\not\in E(G)$. If $xv_3\not\in E(G)$ then $x$ and $v_2v_3$ form an induced subgraph $K_1\cup P_2$, contradicting the hypothesis. Hence $xv_3\in E(G)$. By a similar argument, $xv_i\in E(G)$ for each $i=1,3,5,...t-1$ and $t$ is even. Further, since $C$ is extreme, $\{x,v_2,v_4,...,v_t\}$ is an independent set of vertices. Moreover, for each $u,v\in \{x,v_2,v_4,...,v_t\}$, there is no a path connecting $u$ and $v$ and passing through $V(G)\backslash (V(C)\cup \{x\})$. If $xz\in E(G)$  for some $z\in V(G)\backslash (V(C)\cup \{x\})$ then recalling that $C$ is extreme, we conclude that $zv_2\not\in E(G)$ and $zv_3\not\in E(G)$, contradicting the fact that $G$ is $K_1\cup P_2$-free. Hence  $V(G\backslash C)$ is an independent set of vertices. If $zv_2\in E(G)$ then as above, $zv_4\in E(G)$ and hence 
$$
v_1xv_3v_2zv_4\overrightarrow{C}v_1
$$
is longer than $C$, a contradiction. This means that for each $v\in V(G\backslash C)$, we have $N(v)=\{v_1,v_3,...,v_{t-1}\}$. Put
$$
V_1=V(G\backslash C) \cup \{v_2,v_4,...,v_t\}, \ \ V_2=\{v_1,v_3,...,v_{t-1}\}.
$$
Since $V_1$ is independent and $|V_1|>|V_2|$, $G$ is not hamiltonian and $\alpha(G)>n/2$. If $V_2$ is independent or $G[V_2]$ is complete then $G\in \aleph$ and we are done. Otherwise denote by $V_3$ a largest independent subset in $V_2$ and put $V_4=V_2\backslash V_3$. Let $w_1\in V_3$ and $w_2\in V_4$. Clearly $w_2w_3\in E(G)$ for some $w_3\in V_3\backslash \{w_1\}$, since otherwise $V_3\cup \{w_2\}$ is an independent set of vertices, contradicting the maximality of $V_3$. If $w_1w_2\not\in E(G)$ then $w_1$ and $w_2w_3$ form an induced subgraph $K_1\cup P_2$, contradicting the hypothesis. Hence $w_1w_2\in E(G)$ implying that $N(v)=V_4$ for each $v\in V_3$. Applying the same arguments to $V_4$ instead of $V_2$, we conclude that $G\in\aleph$.   \  \  \  \rule{7pt}{6pt}                    \\

\noindent \textbf{Proof of Theorem 2}. Let $G$ be a 1-tough $(K_1\cup P_3)$-free graph. Since $G$ is 1-tough, it contains a cycle. Let $C$ be a longest cycle in $G$ and $H$ a connected component of $G\backslash C$ of maximum order. If $V(H)=\emptyset$ then $C$ is a Hamilton cycle and we are done. Let $V(H)\not=\emptyset$ and let $\xi_1,...,\xi_s$ be the elements of $N_C(H)$ occuring on $C$ in a consecutive order. Since $G$ is 1-tough, we have $s\ge2$. Set 
$$
I_i=\xi_i\overrightarrow{C}\xi_{i+1},       \     I_i^\ast=\xi_i^+\overrightarrow{C}\xi_{i+1}^-   \   \  (i=1,2,...,s),
$$
where $\xi_{s+1}=\xi_1$. The segments  $I_1,I_2,...,I_s$ are called elementary segments on $C$ induced by $N_C(H)$. We call a path $L=z\overrightarrow{L}w$ an intermediate path between two distinct elementary segments $I_a$ and $I_b$ if
$$
z\in V(I_a^\ast),   \    w\in V(I_b^\ast),    \    V(L)\cap V(C\cup H)=\{z,w\}.
$$

Define $\Upsilon(I_{i_1},I_{i_2},...,I_{i_t})$ to be the set of 
all intermediate paths between elementary segments 
 $I_{i_1},I_{i_2},...,I_{i_t}$. If $\Upsilon(I_1,I_2,...,I_s)=\emptyset$ 
then $G\backslash \{\xi_1,...,\xi_s\}$ has at least $s+1$
 connected components, contradicting the fact that $G$ is 1-tough. 
Otherwise $\Upsilon(I_a,I_b)\not=\emptyset$ for some distinct $a,b\in \{1,...,s\}$. Choose a path $L=x\overrightarrow{L}y$ in $\Upsilon(I_a,I_b)$ such that $x\in V(I_a^\ast)$ and  $y\in V(I_b^\ast)$. If $|V(L)|\ge 3$ then each vertex $v\in V(H)$ with $x\overrightarrow{L}x^{++}$ forms an induced subgraph $K_1\cup P_3$, contradicting the hypothesis. Let $|V(L)|=2$, i.e. $L=xy$. Put 
$$
Q=\xi_a^+\overrightarrow{C}xy\overleftarrow{C}\xi_b^+.
$$
Assume without loss of generality that $L$ is chosen from $\Upsilon(I_a,I_b)$ such that $|V(Q)|$
is minimum. Since $C$ is extreme, by standard arguments, $\{\xi_1,...,\xi_s\}^+$ is an independent set of vertices, implying that either $x\not =\xi_a^+$ or $y\not=\xi_b^+$, say $x\not =\xi_a^+$. Since $|V(Q)|$ is minimum, we have $x^-y\not\in E(G)$, that is $x^-xy$ forms an induced subgraph $P_3$. But then each vertex $v\in V(H)$ with $x^-xy$ forms an induced $K_1\cup P_3$, again contradicting the hypothesis. Theorem 2 is proved. \  \  \  \rule{7pt}{6pt}

\noindent Institute for Informatics and Automation Problems\\ National Academy of Sciences\\
P. Sevak 1, Yerevan 0014, Armenia\\ E-mail: zhora@ipia.sci.am

\end{document}